\documentstyle[amscd,amssymb,epsf]{amsart}
\newtheorem{thm}{Theorem}[section]
\newtheorem{lem}[thm]{Lemma}
\newtheorem{conj}[thm]{Conjecture}

\newtheorem{cor}[thm]{Corollary}

\theoremstyle{definition}
\newtheorem{defn}[thm]{Definition}

\newtheorem{rem}[thm]{Remark}
\newtheorem{nota}[thm]{Notation}

\makeatletter
\def\theequation{\@arabic\c@equation}
\makeatother
\numberwithin{equation}{thm}

\def\Sym{{\text{\rm Sym}}}
\def\cc{\cal C}

\def\R{{\Bbb R}}

\def\Q{{\Bbb Q}}
\def\Z{{\Bbb Z}}

\def\P{{\Bbb P}}
\def\C{{\Bbb C}}

\def\whsq{\vbox to 5.8pt
{\offinterlineskip\hrule
\hbox to 5.8pt{\vrule height
5.1pt\hss\vrule height 5.1pt}\hrule}}

\def\maps{\longrightarrow}

\def\C{{\Bbb C}}
\def\R{{\Bbb R}}
\def\P{{\Bbb P}}

\def\Z{{\Bbb Z}}
\def\Q{{\Bbb Q}}

\def\({\left(}
\def\){\right)}

\def\part{P(n)}

\def\<{\langle}
\def\>{\rangle}

\begin{document}

\title[Modular forms and Donaldson invariants]
{Modular forms and Donaldson invariants for 4-manifolds with $b_+=1$}
\keywords{Blowup-formulas, Donaldson invariants, modular forms}
\author{Lothar G\"ottsche}
\address{Max--Planck--Institut f\"ur Mathematik\\Gottfried--Claren--Stra\ss e
26\\
D-53225 Bonn, Germany}
\email{lothar@@mpim-bonn.mpg.de}

\maketitle\

%\begin{abstract}
%\end{abstract}

\section{Introduction}
The  Donaldson invariants
of a smooth simply connected $4$-manifold $X$ depend by definition on the
choice
of a Riemannian metric $g$.
In case $b^+(X)>1$ they  turn out to be  independent of the metric
as long as it is generic, and thus  give $C^\infty$-invariants of $X$.

We study the case
 $b_+(X)=1$, where the invariants have been introduced  in [Ko].
We denote by $\Phi^{X,g}_{c_1,N}$
the  Donaldson invariant of $X$ with respect to a
lift $c_1\in H^2(X,\Z)$ of $w_2(P)$ for an $SO(3)$ bundle $P$ on $X$
with $-p_1(P)-3=N$.
Kotschick and Morgan showed in \cite{K-M}  that
the invariants only depend on the chamber of the period point of $g$
in
the positive cone $H^2(X,\R)^+$ in $H^2(X,\R)$.
For two metrics $g_1,g_2$, which do not lie on a wall
they express $\Phi^{X,g_1}_{c_1,N}-\Phi^{X,g_2}_{c_1,N}$
as the sum over certain wall-crossing terms $\delta^X_{\xi,N}$,
where $\xi$ runs over all classes in $H^2(X,\Z)$ which define a wall
between $g_1$ and $g_2$.
They also make the following
conjecture.

\begin{conj}\label{KMconj}
\cite{K-M} $\delta^X_{\xi,N}$ is a polynomial in the
multiplication by $\xi$ and the quadratic form $Q_X$
on $H_2(X,\Z)$ whose coefficients depend only on
$\xi^2$, $N$ and the homotopy type of $X$.
\end{conj}
 John Morgan and Peter Ozsv\'ath have told me that they are now able to
prove the  conjecture
\cite{M-O}.

In previous joint papers \cite{E-G1},\cite{E-G2} with
Geir Ellingsrud  we have studied the wall-crossing terms $\delta^S_{\xi,N}$
in the case of algebraic surfaces $S$  with $p_g=0$.
In  \cite{E-G1} we expressed (for so called good walls)
the $\delta^S_{\xi,N}$ in terms of Chern classes of some "standard"
bundles on Hilbert schemes of points on $S$, and proceeded to compute
the leading $6$ terms of $\delta^S_{\xi,N}$
(similar results were also  obtained in \cite{F-Q}).
In \cite{H-P} a Feynman path integral aproach  to this problem is developed,
and some of the leading terms of the wall-crossing formulas are determined.

In \cite{E-G2}, which builds on \cite{E-G1}, we restrict to
the case of rational surfaces and use the Bott residue formula to
compute the $\delta^S_{\xi,N}$ explicitly (with help of the computer).
As an application, using also the blowup formulas,
 we computed e.g. the Donaldson invariants of
$\P_2$ of degree smaller then $50$.

In \cite{K-L} the wall-crossing formulas had already been used in
combination with the blowup formulas to compute  Donaldson invariants
of $\P_2$ and $\P_1\times\P_1$ and to show in particular that neither $\P_2$
nor $\P_1\times \P_1$ is of simple type. Their calculations also showed that
the blowup formulas impose  restrictions on the wall-crossing formulas,
although this is not pursued systematically there. The authors
did however expect that this can be used to determine many
(and possibly all) the wall-crossing formulas for rational surfaces.

In the current paper we  want to show that in fact, assuming
conjecture \ref{KMconj}, one can determine the $\delta^X_{\xi,N}$
completely for all  $X$
and all walls in $H^2(X,\R)^+$ by use of the blowup formulas. We will
 determine
  a universal generating function  $\Lambda(L,Q,x,t,\tau)$
 which expresses all  $\delta^X_{\xi,N}$ for all $X$, $N$ and $\xi$.
Here $\tau$ is a parameter from the complex upper half plane,
$L$,$Q$ and $x$ stand for the multiplication with
$\xi$, the quadratic form and the class of a point, and the exponent
of $t$ is the signature of $X$.
It turns out that $\Lambda(L,Q,x,t,\tau)$ is an exponentional expression
in  certain modular forms (with respect to $\tau$).
As an application of our results we also get a modular forms expressions
for all the Donaldson invariants of the projective plane $\P_2$.
Already in \cite{K-L} it had been shown that the Donaldson invariants
of $\P_2$ and $\P_1\times \P_1$ are determined by the wall-crossing formulas
on the blowup of $\P_2$ in two points. We use instead
a simple fact due to Qin:  on a rational ruled surface the Donaldson
invariants  with respect to a  first Chern class $c_1$ with odd restriction
to a fibre vanish for a special chamber $\cc_F$.

The results of this paper
should be seen in comparison with the new developments of
Seiberg-Witten theory \cite{S-W},\cite{W1} which suggest  a connection between
the Donaldson invariants
and modular forms: the Donaldson invariants (and also the
Seiberg-Witten invariants) are seen as  degenerations of supersymmetric
theories,  parametrized by the  "$u$-plane" (i.e. the modular curve
$\Bbb H/\Gamma(2)$). In fact Witten  informed me that he is currently
trying to determine wall-crossing formulas and the Donaldson invariants
of the projective plane by integrating over the $u$-plane
(see also \cite{W2}). The results should also be related to
the current work  \cite{P-T} towards proving the conjectural relationship
between Seiberg-Witten and Donaldson invariants.

The main tool for getting our result are the blowup formulas,
which   for $4$-manifolds with $b_+=1$ I learned from \cite{K-L}.
Let $\widehat X:=X\#\overline\P_2$
(e.g. if $X$ is an algebraic surface, we can take $\widehat X$
to be the blowup of $X$ at a point).
The idea is very simple:
If $\cc$ is a chamber in $H^2(X,\R)^+$
and $\widehat \cc$ is a related chamber (see below), then there is
a formula relating
the Donaldson invariants of $X$ with respect to $\cc$
and those of $\widehat X$ with respect to $\widehat \cc$.
So let now $\cc_-$ and $\cc_+$ be two chambers separated by the wall
$W^\xi$, then in general there are several walls between the related chambers
$\widehat \cc_-$
and $\widehat \cc_+$ on $\widehat X$, but it is very easy to determine them.
We can therefore express the wall-crossing term   $\delta^X_{\xi,N}$
as follows.
We apply the blowup formulas
to the related chambers $\cc_-,\widehat \cc_-$ and
$\cc_+,\widehat \cc_+$
and add up the wall-crossing terms for all  walls between
$\widehat\cc_-$ and $\widehat \cc_+$.
This gives recursive relations. After encoding our information
into a generating function $\Lambda_X(L,Q,x,t,\tau)$, these recursive
relations translate into differential equations,
which enable us  to determine  $\Lambda_X$ up to multiplication by
a universal function $\lambda(\tau)$.
Unlike the case of the blowup formulas in \cite{F-S},
the modular forms enter the formulas already as the coefficients
of the differential equations; they arize as theta functions
for lattices describing the walls between related chambers.

In order to finally determine $\lambda(\tau)$ we
consider the  particular case $X=\P_1\times\P_1$.
The above mentioned result of Qin now says that,  for first Chern class
$c_1=F+G$,
the sum of the classes of the fibres in the two different directions,
there are always two different chambers $\cc_F$, $\cc_G$ of type $(c_1,N)$
where the corresponding Donaldson invariants vanish.
Therefore the sum of the $\delta^X_{\xi,N}$ for all classes $\xi$ defining
walls between $\cc_F$ and $\cc_G$ must be zero.
This fact gives us an additional
recursion relation, and with this we can finally determine $\lambda(\tau)$.

If we assume only a weaker form of the conjecture, namely if we allow
$\delta^X_{\xi,N}$ to depend on $X$, rather than just on the
homotopy type, then we still get our result for  $X$ a rational surface.
If  we assume the conjecture and the blowup formulas
also in the
case that $X$ is not simply connected but $b_1(X)=0$, then we  can
partially extend our result also to  this case.

I am very thankful to  Don Zagier, who
proved lemma \ref{residue} for me.
I would also like to  thank John Morgan and Stefan Bauer
for very useful conversations.
I would like to thank Dieter Kotschick for sending me the preprint
\cite{K-L}, which was very important both for \cite{E-G2} and for this work.
This paper grew out of the joint work \cite{E-G1}, \cite{E-G2}
with Geir Ellingsrud. Motivated by this work, and based also on
\cite{K-L}, I slowly realized the importace of the blowup formulas
in this context. Also the explicit formulas for the wall-crossing in
\cite{E-G2} were very important for me
to  keep  confidence in my computations.

\section{Background material}

In this paper we will denote by $X$ a simply connected smooth $4$-manifold with
$b_+(X)=1$ and $b_2(X)\ge 2$. We will assume conjecture \ref {KMconj}.
\begin{nota}
For elements $A\in H^2(X,\Q)$ and $\alpha\in H_2(X,\Q)$
we denote by $A\cdot \alpha\in \Q$ the canonical pairing,
by $\check A\in H_2(X,\Z)$ the Poincar\'e dual
and by $A^2$ the number
$A\cdot\check A$.
We denote by $Q_X$ the quadratic form on $H_2(X,\Z)$ and, for a
class $\eta\in H^2(X,\Q)$, by $L_\eta$ the  linear form $\alpha\mapsto
\eta\cdot\alpha$ on $H_2(X,\Q)$.
If there is no risk of confusion we denote by $a$ the reduction of $A\in
H^2(X,\Z)$ modulo $2$.

For a smooth four-manifold $X$ we denote by $\widehat X$
the connected sum $X\#\overline\P_2$ of $X$ with $\P_2$ with the reversed
orientation, (e.g. if $X$ is a smooth complex surface, then
$\widehat X$ is the blowup of $X$ in a point).
Let $E$ be the image of the  generator of $H^2(\overline\P_2,\Z)$
in $H^2(\widehat X,\Z)$. We will will always identify
$H_2(X,\Z)$ with the kernel of $L_E$ on
$H_2(\widehat X,\Z)$. We  write $e$ for the reduction of $E$ modulo $2$.

 Let $g$ be a Riemannian metric on $X$, and $P$
and $SO(3)$ principal bundle with first Pontrjagin class
$p_1(P)=-(N+3)$.
We denote by $\Phi^{X,g}_{c_1,N}$ the Donaldson invariant
corresponding to $P$, the metric
$g$ and the lift $c_1\in H^2(X,\Z)$ of $w_2(P)$.
We use the conventions of e.g. \cite{F-S} which coincide up
to a power of $2$ with the conventions of \cite{Ko}.
If $X$ is an algebraic surface and $H$ an ample divisor
we will write $\Phi^{X,H}_{c_1,N}$ for the invariant with respect to the
Fubini-Studi metric induced by $H$.
Let $p\in H_0(X,\Z)$ be the class of a point.
Let $A_N(X)$ be the set of polynomials of weight $N$ in $H_2(X,\Q)\oplus
H_0(X,\Q)$, where $\alpha\in H_2(X,\Q)$ has weight $1$ and $p$ has weight
$2$. Then $\Phi^{X,g}_{c_1,N}$ is a linear map $A_N(X)\maps \Q$.
We put $\Phi^{X,g}_{c_1,N}:=0$ if $N$ is not congruent to $-c_1^2+3$ modulo $4$
and
$\Phi^{X,g}_{c_1}:=\sum_{N\ge 0}\Phi^{X,g}_{c_1,N}$.
\end{nota}

\subsection{Walls and chambers}

\begin{defn}\label{defwall} (see e.g. \cite{Ko}, \cite{K-M})
Let $w\in H^2(X,\Z/2\Z)$
and  $N$  a nonnegative integer.
 Let $H^2(X,\R)^+$ be the positive  cone in $H^2(X, \R)$.
For $\xi\in H^2(X,\Z)$ let
$$W^\xi:= \{ x\in H^2(X, \R)^+ \bigm| \xi\cdot\check x =0\big\}.$$
We shall call $W^\xi$ a wall of type $(w,N)$, and say that it is
defined by  $\xi$, if  $w$ is the reduction  of $\xi$ modulo $2$,
$N+3$ is congruent to $\xi^2$ modulo $4$
 and $-(N+3)\le \xi^2<0$.
Note that any class $\xi\in H^2(X,\Z)$ with $\xi^2<0$
will define a wall of type $(w,N)$ for suitable $N$ and $w$ the reduction
of $\xi$ modulo $2$; we will in this case say that $\xi$
defines a wall of type $(N)$.
A chamber of type $(w,N)$ is a connected component of
the complement of the walls of type $(w,N)$ in $H^2(X, \R)^+$.

For a  Riemannian metric $g$ on $X$ we denote by $\omega(g)\in H^2(X,\R)^+$
the corresponding period point.
If $(w,N)$ are given, a metric is called generic if its period point
$\omega(g)$
does not lie on a wall of type $(w,N)$.
For $A_-,A_+\in H^2(X,\R)$  we denote
by
$W_{w,N}^X(A_-,A_+)$ the set of all $\xi\in H^2(X,\Z)$ defining a wall
of type $(w,N)$ with
$\xi\cdot \check A_-<0<\xi\cdot \check A_+$.
We put
$$W_{w}^X(A_-,A_+):=\bigcup_{N\ge 0} W_{w,N}^X(A_-,A_+).$$
\end{defn}

\begin{thm} \label{wallchange}\cite{K-M}
Let $c_1\in H^2(X,\Z)$ and $w$ the reduction of $c_1$ modulo $2$.
For all $\xi\in H^2(X,\Z)$ defining a wall of type $(w,N)$
we put $\varepsilon(c_1,\xi,N):=(5N+3+\xi^2+(\xi-c_1)^2)/4$.
There exists
$\delta^X_{\xi,N}:\Sym^N(H_2(X,\Q))\maps \Q$ such that for all generic metrics
$g_+$ and $g_-$ with  $\omega(g_+)$ and $\omega(g_-)$ in the same connected
component of $H^2(X,\R)^+$
$$\Phi^{X,g_+}_{c_1,N}-\Phi^{X,g_-}_{c_1,N}=
\sum_{\xi\in W_{w,N}^X(\omega(g_-),\omega(g_+))}
(-1)^{\varepsilon(c_1,\xi,N)}\delta^X_{\xi,N}.$$
Furthermore, if $\omega(g_1)=-\omega(g)$, then
$\Phi^{X,g_1}_{c_1,N}=-\Phi^{X,g}_{c_1,N}$.
\end{thm}

\begin{rem}\label{sign}
\begin{enumerate}
\item Our sign conventions are different from those of \cite{K-M}
and \cite{K-L}. In fact the sign is chosen in order to give
the leading term $L_\xi^{N-2d}Q_X^d$ (with $d=(N+3+\xi^2)/4$)
of $\delta^X_{\xi,N}$
a positive coefficient.
\item In the future we will always  implicitely assume that all the metrics
that
we consider have their period point in the same connected component of
$H^2(X,\R)^+$.
\item By theorem \ref{wallchange} we can write $\Phi^{X,\cc}_{c_1,N}:=
\Phi^{X,g}_{c_1,N}$ for any metric $g$ with $\omega(g)$ in the chamber
$\cc$.
\end{enumerate}
\end{rem}

\subsection{Blowup formulas}
The blowup formulas relate the Donaldson invariants of a $4$-manifold $Y$
and $\widehat Y=Y\#\overline\P_2$.
In the case $b_+(Y)>1$, when the invariants do not depend on the chamber
structure, they have been shown e.g. in \cite{O}, \cite{L}
and  in the most general form in \cite{F-S}. In the case when
$X$ is a simply
connected $4$-manifold with $b_+=1$
I learned the blowup formulas from \cite{K-L}.
They then  depend on the chamber structure.

\begin{defn}\label{related}(see \cite{Ko}).
Let $\cc\subset H^2(X,\R)^+$ be a chamber of type $(w,N)$.
A chamber $\cc_0\subset H^2(\widehat X,\R)^+$ of type $(w,N)$ (resp.
$\cc_e\subset H^2(\widehat X,\R)^+$ of type $(w+e,N+1)$ is said to be related
to $\cc$ if and only if
$\cc$ is contained in the closure
$\overline \cc_0$ (resp in $\overline\cc_e$).
\end{defn}

By \cite{T} the formulas of \cite{F-S} also hold for $X$ with $b_+(X)=1$,
we will however only need a quite easy special case (see e.g. \cite{Ko}
and \cite{S}, ex. 11).

\begin{thm} \label{blowup}
Let  $\cc\subset H^2(X,\R)^+$ be a chamber of type $(w,N)$, and
let  $\cc_0\subset H^2(\widehat X,\R)^+$  (resp.
$\cc_e\subset H^2(\widehat X,\R)^+$) be related chambers of types $(w,N)$
(resp. $(w+e,N+1)$).
Then  for all $\alpha\in A_N(X)$ and $\beta\in A_{N-2}(X)$
for which both sides are defined we have
\begin{align}\tag*{$(0)_b$}
 \Phi^{X,\cc}_{c_1,N}(\alpha)&=\Phi^{\widehat X,\cc_0}_{c_1,N}(\alpha),\\
\tag*{$(1)_b$}\Phi^{X,\cc}_{c_1,N}(\alpha)
&=\Phi^{\widehat X,\cc_e}_{c_1+E,N+1}(\check E\alpha),\\
\tag*{$(2)_b$}\Phi^{\widehat X,\cc_0}_{c_1,N}(\check E^2\beta)&=0,\\
\tag*{$(3)_b$}\Phi^{X,\cc}_{c_1,N}(x\beta)&=-\Phi^{\widehat X,
\cc_e}_{c_1+E,N+1}(\check E^3\beta).
\end{align}
\end{thm}

\subsection{Extension of wall-crossing formulas}
We want to extend theorem \ref{wallchange}
from $\Sym^N(H_2(X,\Q))$ to $A_N(X)$.
For this we have to extend the definition of $\delta^X_{\xi,N}$.
In the case that $\xi$ is divisible by $2$
(i.e.  $w=0$)
  we also have to extend the definition
of $\Phi^{X,g}_{c_1,N}$  to classes not in the stable range.
For technical reasons we also redefine the
$\delta^X_{\xi,N}$ in case the intersection form on $H_2(X,\Z)$ is even
or the rank of $H_2(X,\Z)$ is at most $2$. It should not be difficult
to prove that this definition agrees with that of \cite{K-M},
but we only need that theorem \ref{wallchange}
still holds.

\begin{defn}\label{extend}
\begin{enumerate}
\item Let
 $N=4c_2-3$ for $c_2\in \Z$.
Let $\cc$ be a chamber of type $(0,N)$ in $H^2(\widehat X,\R)^+$
and $\cc_e$ a related chamber of type $(e,N+1)$ on $\widehat X$.
Then we put for all $\alpha\in A_N(X)$
$$\Phi^{X,\cc}_{0,N}(\alpha):=\Phi^{\widehat X,\cc_e}_{E,N+1}(\check
E\alpha).$$
Note that $(1)_b$ above guaranties that our definition
restricts to the
standard definition if $\alpha$ is in the stable range.

\item Let $\xi\in H^2(X,\Z)$ with $\xi^2<0$. We extend the definition of
$\delta^X_{\xi,N}$ by putting $\delta^X_{\xi,N}:=0$ if
$\xi$ does not define a wall of type $(N)$
(i.e. if $\xi^2$ is not congruent ot $N+3$ modulo $4$ or $N+3+\xi^2<0$).

\item
Assume now that the intersection form on $H_2(X,\Z)$ is even
or the rank of $H_2(X,\Z)$ is at most $2$ or that $\xi$ is divisible by
$2$ in $H^2(X,\Z)$. Then we put for $\alpha\in \Sym^N(H_2(X,\Q))$
$$\delta^X_{\xi,N}(\alpha):=\sum_{n\in\Z}(-1)^{n-1}
\delta^{\widehat X}_{\xi+(2n+1)E,N+1}(\check E\alpha).$$
Note that by (1)
the sum runs in fact only through integers $n$ with
$(2n+1)^2\le N+4+\xi^2$.

\item Assume that $\delta^Y_{\eta,N}(p^r\beta)$
is already defined for all $m$ for
 $Y=X\#m\overline \P_2$ for all $N$, all $\xi\in H^2(Y,\Z)$
with $\xi^2<0$ and all
$\beta\in \Sym^{N-2r}(H_2(Y,\Q))$.
Then we put
$$\delta_{\xi,N}^Y(p^{r+1}\alpha):=\sum_{{n\in\Z}}(-1)^{n}
\delta^{\widehat Y}_{\xi+(2n+1)E,N+1}(\check E^3p^r\alpha)$$
for all $\alpha\in \Sym^{N-2r-2}(H_2(Y,\Q))$. Again by (1)
the sum runs  only through $n$ with
$(2n+1)^2\le N+4+\xi^2$.
\end{enumerate}

We note that by definition $\delta_{\xi,N}^X=0$
 if $\xi$ does not define a wall of type
$(N)$. Finally we put
$$\delta^X_\xi:=\sum_{N\ge 0} \delta^X_{\xi,N}.$$
If  $t$ is an indeterminate,
we write
$\delta^X_\xi(\sum_N \alpha_N t^N)$ for $\sum_N \delta^X_{\xi,N}(\alpha_N)
t^N$,
and similarly for $\Phi^{X,g}_{c_1}$.
\end{defn}

\begin{rem}\label{c10}
There is a small subtlety about the definition of $\Phi^{X,\cc}_{0,N}$
in (1). If  $w\ne 0$, then,
given a chamber $\cc$ in $H^2(X,\R)^+$ of type
$(w,N)$, there is a unique related chamber $\cc_e$ in
$H^2(X,\R)^+ $ of type
$(w+e,N+1)$ consisting of all
$ \mu+a E $ with $\mu\in \cc$ and $a\in \R$ sufficiently small.
If $w=0$, however, $E$ defines a wall of type $(e,N+1)$ separating
two chambers $\cc_e^+$ (corresponding to $a>0$) and $\cc_e^-$
(corresponding to $a<0$) of type $(e,N+1)$,
which are both related to $\cc$.
$\Phi_{0,N}^{X,\cc}$ is still well-defined, as
$\delta_{(2n+1)E,N+1}^{\widehat X}(\check E^{2k+1}\alpha)=0$
for $N$ congruent to $1$ modulo $4$: By
conjecture \ref{KMconj}
(and the extension  \ref{point} to $A_N(X)$)
$\delta_{E,N+1}^{\widehat X}(p^r \bullet)$ is a polynomial in $L_E$ and
$Q_{\widehat X}$,  $N+1$ is even and $E\cdot \alpha=0$.

Similarly, if $w\ne 0$, then
there is a unique related chamber $\cc_0$ in
$H^2(X,\R)^+ $ of type
$(w,N)$. If $w=0$, then there are two
related chambers separated by a wall defined by $2E$,
but $\delta_{2nE,N}^{\widehat X}(\check E^{2k}\alpha)=0$.
\end{rem}

\subsection{Vanishing on rational ruled surfaces}
An important rule both in the proof of the main theorem and
in the application to Donaldson invariants of the projective plane
is played by the following elementary vanishing result.
Let $S$ be a rational ruled surface, and let $F,E\in H^2(S,\Z)$ be the classes
of a fibre of the projection to $\P_1$ and a section respectively.
For an ample divisor $H$ let $M^S_{H}(c_1,c_2)$
be the moduli space of $H$-stable torsion-free sheaves with Chern classes
$(c_1,c_2)$.

\begin{lem} \label{qin}\cite{Q2}
Assume  $c_1\cdot F=1$,
then $M^S_{F+\epsilon E}(c_1,c_2)$ is empty for all
sufficiently small $\epsilon >0$. In particular,
given $N\ge 0$, we get $\Phi^{S,F+\epsilon E}_{c_1,N}=0$ for all
sufficiently small $\epsilon >0$.
\end{lem}

\section{Main Theorem}
We want to express the wall-crossing formulas in terms of the $q$-development
of certain modular forms. We start by reminding the reader of
some notations and elementary facts (see e.g. \cite{H-B-J},
\cite{R}).

\begin{nota}
Let $\Bbb H=\{\tau\in \C \,|\, Im(\tau)>0\}$ be the complex
upper half plane. We denote
$q=e^{2\pi i\tau}$ and $q^{1/n}=e^{2\pi i\tau/n}$.
For a positive integer $n$ let
$$\sigma_k(n):=\sum_{d|n}d^k\quad \hbox{and} \quad
\sigma_1^{odd}(n):=\sum_{d|n,\ d\ odd}d^k.$$
Let $\eta(\tau):=q^{1/24}\prod_{n>0} (1-q^n)$ be the Dirichlet eta-function,
and let $\Delta(\tau)=\eta(\tau)^{24}$ be the discriminant.
We denote
$$\theta(\tau):=\sum_{n\in \Z}q^{n^2}$$ the theta function for the latice
$\Z$. We also have the Eisenstein series
$$G_{2}(\tau):=-1/24+\sum_{n\ge 1}\sigma_{1}(n)q^n$$
and the $2$-division value
$$e_3(\tau):=1/12+2\sum_{n\ge 1} \sigma_1^{odd}(n)q^{n/2}$$
We put $f(\tau):=\eta(2\tau)^3/\theta(\tau).$ Then
$\eta(2\tau)$, $\theta(\tau)$, $G_{2}(2\tau)$, $e_3(2\tau)$ and  $f(\tau)$
are modular forms of weights $1/2$, $1/2$, $2k$, $2$ and $1$
respectively for certain subgroups of $SL(2,\Z)$.

We will denote
$ d\log_q(g):=g^{-1}{d g/ dq}$.
Note that
$$\leqno{(*)} \quad  d\log_q(g_1g_2)= d\log_q(g_1)+ d\log_q(g_2)
\hbox{\ \ \ and \ \ \ }
 d\log_q(g_1/g_2)= d\log_q(g_1)- d\log_q(g_2).$$
\end{nota}

\begin{rem}\label{modular} We will  use the following identities
\begin{eqnarray*}
&{(1)\ }&\eta(2\tau)^3=\sum_{n\in\Z} (-1)^n (n+1/2)q^{(n+1/2)^2},\\
&{(2)\ }&\theta(\tau)={\eta(2\tau)^5\over
\eta(\tau)^2\eta(4\tau)^2},\quad\quad\quad\quad\quad\quad
{(3)\ } f(\tau)={\eta(\tau)^2\eta(4\tau)^2\over \eta(2\tau)^2 },\\
&{(4)\ }& q\, d\log_q(\eta(2\tau))=-2G_2(2\tau),\quad\quad\quad
{(5)\ } q\, d\log_q(\theta(\tau))=-2G_2(2\tau)-e_3(2\tau).
\end{eqnarray*}
\end{rem}
\begin{pf} (1) and (2) are standard facts, following e.g. from the Jacobi
identity. (3) follows from (2). (4)  follows by an easy
 calculation using $(*)$, and, using also (2), the proof of (5) is similar.
\end{pf}

The main result of this paper is the following.

\begin{thm} \label{mainthm}
Let $X$ be a simply connected $4$-manifold with $b_+=1$ and signature
$\sigma(X)$.
Let $\xi\in H^2(X,\Z)$ with $\xi^2<0$.
For $\alpha\in H_2(X,\Z)$ denote
\begin{eqnarray*}
g^X_{\xi}(\alpha z,x,\tau)&:=&\exp\Big(({\xi/2}\cdot\alpha)z/
f(\tau)-(Q_X(\alpha)/2)z^2
(2G_2(2\tau)+e_3(2\tau))/f(\tau)^2\\
&&\quad -3xe_3(2\tau)/f(\tau)^2\Big)
\theta(\tau)^{\sigma(X)}f(\tau){\Delta(2\tau)^2
\over\Delta(\tau)\Delta(4\tau)}.\end{eqnarray*}
Then
$$ \delta^X_\xi\big(\exp(\alpha z+px)\big)
=\hbox{res}_{q=0}(q^{-\xi^2/4}g^X_{\xi}(\alpha z,x,\tau){dq/ q}).
$$
\end{thm}

\begin{rem} \label{power}
\begin{enumerate}
\item One can see that this expression for $\delta_{\xi}^X$
is not  compatible with the simple type condition. In particular
given $c\in H^2(X,\Z)$
a $4$-manifold $X$ with $b_+=1$ will  be of c-simple type at most
for some special points in the closure $\overline C_X$
of the positive cone of $X$.
It had already been shown in \cite{K-L} that
$\P_2$ is not of simple type and that there is no chamber for which
$\P_1\times \P_1$ is of simple type. It is easy to see
from this that rational algebraic
surfaces $X$
can be of simple type at most for  special points in $\overline C_X$.
\item
The expression  $\hbox{res}_{q=0}(q^{-\xi^2/4}g^X_{\xi}(\alpha z,x,\tau)dq/q)$
is just the coefficient of $q^{\xi^2/4}$ of $g^X_{\xi}(\alpha z,x,\tau)$.
The current formulation
is however more intrinsic. Note also that
$dq/q=2\pi i d\tau$.
\item We see that the coefficient $g_{N-2r,r}$ of $z^{N-2r}x^r$ in
$g^X_{\xi}(\alpha z,x,\tau)$
 is $q^{-(N+3)/4}$ multiplied with a power series in $q$
In particular, if $\xi$ defines a wall of type $(N)$, then
$q^{-\xi^2/4}g_{N-2r,r}$ is a Laurent series in $q$.
If $\xi$ with $\xi^2<0$ does not define a
wall of type $(N)$,
then the constant term of $q^{-\xi^2/4}g_{N-2r,r}$ is zero.
\item It would be interesting to know whether for
classes $\xi\in H^2(X,\Z)$ with $\xi^2\ge 0$ the expression
$\hbox{res}_{q=0}(q^{-\xi^2/4}g^X_{\xi}(\alpha z,x,\tau)dq/q)$
has a geometrical or gauge-theoretical meaning.
\end{enumerate}
\end{rem}

As a reasonably straightforward application of theorem \ref{mainthm}
we can determine all the Donaldson invariants of the projective plane $\P_2$.

\begin{thm}\label{donp2} We denote by $\sqrt{i}$ a primitive
$8$-th root of unity and by $H$ in $H^2(\P_2,\Z)$ the hyperplane class.
Put $$e_n(z,x,\tau):=\exp\big((n/2) \sqrt{i} z  /f(\tau)-i(z^2/ 2)
(2G_2(2\tau)+e_3(2\tau))/f(\tau)^2-3ix e_3(2\tau)/f(\tau)^2\big)
{\Delta(2\tau)^2\over\Delta(\tau)\Delta(4\tau)}.$$
Then
\begin{align*}
\tag*{$(1)$}
&\Phi_{H}^{\P_2}(\exp(\check Hz+px))=
\hbox{res}_{q=0}\Big(\sum_{{n>0 \ \text{odd}}\atop {a>n}\ \text{even}}
(-1)^{(n+1)/2}
q^{(a^2-n^2)/4}e_n(z,x,\tau)f(\tau)\Big){dq/q}.\\
\tag*{$(2)$}
&\Phi_{0}^{\P_2}(\exp(\check Hz+px))=\hbox{res}_{q=0}
\Big(\sum_{{{n>0}\ \text{even} }\atop{{a>n}\ \text{odd}}} (-1)^{(a-1)/2}
q^{(a^2-n^2)/4}{a\over 2\sqrt{i}}e_n(z,x,\tau)\Big){dq/ q}.
\end{align*}
In $(2)$ we have used definition \ref{extend} to define $\Phi_{0}^{\P_2}(\check
H^{N-2r}p^r)$ for $r\ge (N-5)/4$.
One can check that (up to different sign conventions) $(1)$ and $(2)$  agree
with the explicit computations in \cite{K-L} and \cite{E-G2}.
\end{thm}

\begin{pf} (of theorem \ref{donp2} from theorem \ref{mainthm}).
Let $Y$ be the blowup of $\P_2$ in a point, and let $E\in H^2(Y,\Z)$
be the class of the exceptional divisor. Let $F=H-E$ be the class
of a fibre of the ruling $Y\maps \P_1$.
Fix a nonnegative integer $N$.
By lemma \ref{qin} we get for $\epsilon>0$ sufficiently small
$\Phi^{Y,F+\epsilon E}_{H,N}=0=\Phi^{Y,F+\epsilon E}_{E,N+1}$.
On the other hand the chamber of $H-\epsilon E$ is related to the
polarisation $H$ of $\P_2$. Thus we obtain by the blowup formulas
$(0)_b$ and $(1)_b$ that
$\Phi^{\P_2}_{H,N}=\Phi^{Y,H-\epsilon E}_{H,N}$
and
$\Phi^{\P_2}_{0,N}=\Phi^{Y,H-\epsilon E}_{E,N+1}$.
So we get by theorem \ref{wallchange} (and lemma \ref{wallchange1}
below) the formulas
\begin{eqnarray*}
\Phi_{H}^{\P_2}(\exp(\check Hz+px))&=&\sum_{\xi\in W^Y_h(F,H)}
\sqrt{i}^{(\xi^2+3)+(\xi-H)^2}\delta^Y_{\xi}(\exp(-\sqrt{i}\check Hz+i
px)),\\
 \Phi_{0}^{\P_2}(\exp(\check Hz+px))&=&\sum_{\xi\in W^Y_e(F,H)}
\sqrt{i}^{(\xi^2+3)+(\xi-E)^2}
\delta^Y_{\xi}(-\sqrt{i}\check E \exp(-\sqrt{i}\check Hz+ipx)).
\end{eqnarray*}
It is easy to see that
\begin{eqnarray*}
W^Y_h(F,H)&=&\big\{(2n-1)H-2aE\bigm| a\ge n\in\Z_{>0}\big\},\\
W^Y_e(F,H)&=&\big\{ 2nH-(2a-1)E\bigm| a>n\in\Z_{>0}\big\}.
\end{eqnarray*}
For $\xi=nH-aE$ we get
$-\xi^2/4=(a^2-n^2)/4$. Furthermore
$\sqrt{i}^{(\xi^2+3)+(\xi-H)^2}=(-1)^{(n+1)/2}$ if $n$ is odd and $a$ is even
and $\sqrt{i}^{(\xi^2+3)+(\xi-E)^2}
=i^{a+2}$ if $n$ is even and $a$ is odd.
Thus, replacing  $-\sqrt{i}$ by $\sqrt{i}$, (1) follows directly by applying
theorem \ref{mainthm}.
(2)  follows the same way using
that
\begin{eqnarray*}
\delta^Y_{\xi}(-\sqrt{i}\check E \exp(-\sqrt{i}\check Hz+i
px))&=&{d\over dw} \Big(\delta^Y_\xi(\exp(-\sqrt{i}(\check Ew+\check Hz)+i
px))\Big)\Big|_{w=0}\\
&=&\hbox{res}_{q=0}\Big(q^{-\xi^2/4}
{d\over dw}(g^{\widehat\P_2}_{\xi}(-\sqrt{i}(\check Ew+\check
Hz),ix,\tau))\big|_{w=0}\Big).
\end{eqnarray*}
\end{pf}

\begin{rem}
The arguments in  section 6 of \cite{E-G2} show that, using
theorem \ref{mainthm} and the blowup formulas, we can get
explicit generating functions  for all the Donaldson invariants of all rational
surfaces $S$ in all chambers of $H^2(S,\R)^+$.
In \cite{K-L} it had been shown (also using the blowup formulas)
that the wall-crossing terms on $\P_2\#2\overline\P_2$
determine the Donaldson invariants on $\P_2$ and $\P_1\times \P_1$.
\end{rem}

\section{Proof of the main theorem}
We give a brief outline of the argument.
Let $\xi\in H^2(X,\Z)$ define a wall of type $(w,N)$,
and let $\cc_-$ and $\cc_+$ be the two chambers separated by
$W^\xi$. The related chambers $\cc_{-0}$ and $\cc_{+0}$
of type $(w,N)$ (resp. $\cc_{-e}$ and $\cc_{+e}$
of type $(w+e,N+1)$) on $\widehat X$ are now separated by several
walls. We can express
$\delta_{\xi,N}^X$  by first applying the blowup
formula to the pair $\cc_-$,  $\cc_{-0}$ of related chambers,
then summing up the wall-crossing formulas for all walls between
$\cc_{-0}$ and $\cc_{+0}$ and finally applying again the blowup
formula for $\cc_+$,  $\cc_{+0}$
(and similarly for $\cc_{-e}$, $\cc_{+e}$).
The blowup formulas $(0)_b$--$(3)_b$ from \ref{blowup}
will give  relations $(0)_r$--$(3)_r$
between the $\delta_{\xi,N}^X$ and the $\delta_{\xi,N}^{\widehat X}$.
Using conjecture \ref{KMconj}
we  encode this information (for all blowups of $X$)
in a suitable generating function $\Lambda_X$ in several variables.
Then we can
translate $(0)_r$--$(3)_r$
 into differential equations
$(0)_d$--$(3)_d$ for $\Lambda_X$,
which determine $\Lambda_X$ up to multiplication by a
function $\lambda_X(\tau)$.
We finally determine $\lambda_X(\tau)$ by specializing to the case
$X=\P_1\times \P_1$ and applying lemma \ref{qin}.

\begin{lem} \label{wallchange1}  Let $w$ be the reduction modulo $2$ of
$c_1\in H^2(X,\Z)$, and let $N$ be a nonnegative integer.
Let $g_-$ and $g_+$ be two metrics on $X$, whose period points
$\omega(g_-)$ and $\omega(g_+)$  do not lie on a wall
of type $(w,N)$. We denote
$W:=W^X_{w,N}(\omega(g_-),\omega(g_+))$. Then we have for all $\alpha\in
A_N(X)$ and $\beta\in A_{N-2}(X)$:
\begin{align}
\Phi^{X,g_+}_{c_1,N}(\alpha)-\Phi^{X,g_-}_{c_1,N}(\alpha)&
=\sum_{\xi\in W}
(-1)^{\varepsilon(c_1,\xi,N)}\delta^X_{\xi,N}(\alpha)\tag*{$(a)$},\\
\tag*{$(0)_r$} \Phi^{X,g_+}_{c_1,N}(\alpha)-\Phi^{X,g_-}_{c_1,N}(\alpha)&
=\sum_{\xi\in W} (-1)^{\varepsilon(c_1,\xi,N)}\sum_{n\in \Z}\delta^{\widehat
X}_{\xi+2nE}(\alpha),\\
\tag*{$(1)_r$}\Phi^{X,g_+}_{c_1,N}(\alpha)-\Phi^{X,g_-}_{c_1,N}(\alpha)&
=\sum_{\xi\in W} (-1)^{\varepsilon(c_1,\xi,N)}\sum_{n\in
\Z}(-1)^{n-1}\delta^{\widehat X}_{\xi+(2n+1)E,N+1}(\check E\alpha),\\
\tag*{$(2)_r$} 0&=\sum_{\xi\in W} (-1)^{\varepsilon(c_1,\xi,N)}\sum_{n\in
\Z}\delta^{\widehat X}_{\xi+2nE,N}(\check E^2\beta),\\
\tag*{$(3)_r$} \Phi^{X,g_+}_{c_1,N}(p\beta)-\Phi^{X,g_-}_{c_1,N}(p\beta)&
=\sum_{\xi\in W} (-1)^{\varepsilon(c_1,\xi,N)}
\sum_{n\in \Z}(-1)^{n}\delta^{\widehat X}_{\xi+(2n+1)E,N+1}(\check E^3\beta).
\end{align}
(a) says that theorem \ref{wallchange} extends to our
definition of $\delta^X_{\xi,N}$.
\end{lem}

\begin{pf}
We assume that $N$ is congruent to  $3-c_1^2$ modulo $4$
(otherwise both sides of $(a)$, $(0)_r$--$(3)_r$
are trivially zero).
Let $\cc_-$ and $\cc_+$ be the chambers of type $(w,N)$ of
$\omega(g_-)$ and $\omega(g_+)$ respectively.
Let $\cc_{-0}$ and $\cc_{+0}$ (resp. $\cc_{-e}$ and $\cc_{+e}$)
be related chambers in $H^2(\widehat X,\R)^+$ of type $(w,N)$
(resp. $ (w+e,N+1)$).

\noindent{\it Claim:}
\begin{eqnarray*}W_{w,N}^{\widehat  X}(\cc_{-0},\cc_{+0})&=&
\big\{\xi+2nE\bigm |\xi\in W,\ n\in \Z,\ n^2\le (N+3+\xi^2)/4\big\},\\
W_{w+e,N+1}^{\widehat  X}(\cc_{-e},\cc_{+e})&=&
\big\{\xi+(2n+1)E\bigm| \xi\in W,\ n\in \Z,\ (2n+1)^2\le N+4+\xi^2\big\}.
\end{eqnarray*}
In the case that $w=0$, we assume that
$\cc_{-0}$ and $\cc_{+0}$ (resp. $\cc_{-e}$ and $\cc_{+e}$) lie on the same
side of $W^{2E}$ (resp. $W^E$).
The claim is essentially obvious:
Any  $\eta\in W_{w,N}^{\widehat  X}(\cc_{-0},\cc_{+0})$
must be of the form $\xi+\alpha E$ for $\xi\in W$.
By the definition of a wall we see that $\alpha$ must be an even integer
$2n$ with $n^2\le (N+3+\xi^2)/4$. On the other hand it is obvious that
all $\xi+2n E$ with $\xi\in W$ and
$n^2\le (N+3+\xi^2)/4$ lie in $W_{w,N}^{\widehat  X}(\cc_{-0},\cc_{+0})$.
For $W_{w+e,N+1}^{\widehat  X}(\cc_{-e},\cc_{+e})$ we argue analogously.

Using this description of $W_{w,N}^{\widehat  X}(\cc_{-0},\cc_{+0})$ and
$W_{w+e,N+1}^{\widehat  X}(\cc_{-e},\cc_{+e})$,
we see that for $X$ with $b_2(X)>2$ and odd intersection form
and $\alpha\in \Sym^N(H^2(X,\Q))$ and $\beta\in \Sym^{N-2}(H^2(X,\Q))$,
the formulas $(0)_r$--$(3)_r$
are just straightforward translations of
$(0)_b$--$(3)_b$ (note that
$\varepsilon(c_1,\xi,N)-\varepsilon(c_1+E,\xi+(2n+1)E,N+1)$ is congruent to
$n-1$ modulo $2$).

Now assume that  the intersection form of $X$ is even or
$b_2(X)\le 2$ or $\xi$ is divisible by $2$ in $H^2(X,\Z)$. Then definition
\ref{extend}, the description of
$W_{w,N}^{\widehat  X}(\cc_{-e},\cc_{+e})$
and $(1)_b$ imply immediately that $(a)$ holds for all $\alpha\in
\Sym^N(H^2(X,\Q))$.
We  show $(0)_r$--$(3)_r$ for  $\alpha\in \Sym^N(H^2(X,\Q))$
(we only
carry out the case of $(0)_r$, the other cases are analogous.)
Let $\widetilde X:=\widehat X\# \overline \P_2$, we denote
by $F$ the generator of $H^2(\overline \P_2,\Z)$.
Then by definition \ref{extend} and $(0)_r$ for $\widehat X$ we get
\begin{eqnarray*}
\Phi^{X,\cc_+}_{c_1,N}(\alpha)-\Phi^{X,\cc_-}_{c_1,N}(\alpha)
&=&\sum_{\xi\in W} \sum_{m\in \Z}  (-1)^{m+1}
\delta^{\widehat X}_{\xi+(2m+1)E,N+1}(\check E\alpha)\\
&=&\sum_{\xi\in W}\sum_{n\in \Z} \sum_{m\in \Z} (-1)^{m+1}
\delta^{\widetilde X}_{\xi+2nF+(2m+1)E,N+1}(\check E\alpha)\\
&=&\sum_{\xi\in W}
\sum_{n\in \Z} \delta^{\widehat X}_{\xi+2nE,N}(\alpha).
\end{eqnarray*}

Now let $X$  be general. We assume
$(a)$, $(0)_r$--$(3)_r$  for all blowups $Y$ of $X$
and all classes $\alpha=p^l\beta$
with $\beta\in \Sym^{k}(H_2(Y,\Q))$ for some $k$.
Then  $(3)_r$ implies immediately
$(a)$  for $p\alpha$. The proof of
$(0)_r$--$(3)_r$ for $p\alpha$
is analogous to  the last section. We  only carry out the
case of $(1)_r$.
Let $\widetilde Y:=\widehat Y\# \overline \P_2$, we denote
by $F$ the generator of $H^2(\overline \P_2,\Z)$. We get by definition
\ref{extend}
\begin{eqnarray*}
\Phi^{Y,\cc_+}_{c_1,N}(p\alpha)-\Phi^{Y,\cc_-}_{c_1,N}(p\alpha)
&=&\sum_{\xi\in W}
\sum_{n\in \Z}(-1)^{n} \delta^{\widehat Y}_{\xi+(2n+1)E,N+1}
(\check E^3\alpha)\\
&=&\sum_{\xi\in W}\sum_{n\in \Z} \sum_{m\in \Z} (-1)^{n}(-1)^{m-1}
\delta^{\widetilde Y}_{\xi+(2n+1)E+(2m+1)F,N+2}(\check F\check E^3\alpha)\\
&=&\sum_{\xi\in W} \sum_{m\in \Z}  (-1)^{m-1}
\delta^{\widehat Y}_{\xi+(2m+1)E,N+1}(\check Ep\alpha).
\end{eqnarray*}
\end{pf}

\begin{lem}\label{combin}
For $\xi\in H^2(X,\Z)$ we get
$$\exp\big(L_{(\xi+nE)/2}+Q_{\widehat X}\big)(\check E^k\bullet)=
\sum_{s+2t=k} (n/2)^s(-1)^{s+t}{k!\over s!t!}
\exp(L_{\xi/2}+Q_{ X}),$$
as a map $\sum_{N\ge 0} \Sym^N(H_2(X,\Q))\maps \Q$.
\end{lem}
\begin{pf}
$$
\exp\big(L_{(\xi+nE)/2}+Q_{\widehat X}\big)(\check E^k\bullet)=
{d^k\over dw^k}
\exp\big((L_{\xi}-nw)/2+Q_X-w^2\big)\big|_{w=0},$$
and the result follows by induction.
\end{pf}

\begin{rem}\label{point}
Using definition \ref{extend}, lemma \ref{combin}
and easy induction we see that conjecture \ref{KMconj}
implies that
$\delta^X_{\xi,N}(p^r\bullet)$ is a polynomial in $L_{\xi/2}$ and
$Q_X$ with coefficients only depending on $N$, $\xi^2$, $r$
and the homotopy type of $X$.
\end{rem}

\begin{defn}\label{pcoeff}
For all $b\ge 0$ let $X(b):=X\# b\overline\P_2$.
Let $l,k,r,b\in\Z$, put $N:=l+2k+2r$, and assume that there exists a
class $\xi\in H^2(X(b),\Z)$ with $w=\xi^2/4<0$. Then  we put
$$P(l,k,r,b,w):={l!k!\over
(l+2k)!}\hbox{Coeff}_{L_{\xi/2}^lQ_{X(b)}^k}\delta^{X(b)}_{\xi,N}(p^r\bullet).
$$
(By definition $P(l,k,r,b,w)$ will be zero if $\xi$ does not define a wall
of type $(N)$ or if one of $l,k,r,b$ is negative).
Note that $P(l,k,r,b,w)$ is well defined:
By conjecture \ref{KMconj} and remark \ref{point}
$\delta^{X(b)}_{\xi,N}(p^r\bullet)$ is a polynomial in
$L_{\xi/2}$ and $Q_{X(b)}$. As $b_2(X)>1$, the monomials
$L_{\xi/2}^lQ_{X(b)}^k$ are linearly independent as linear maps
$\Sym^{l+2k}(H_2(X,\Q))\maps \Q$, therefore the coefficients of
$L_{\xi/2}^lQ_{X(b)}^k$ in $\delta^{X(b)}_{\xi,N}(p^r\bullet)$
are well-defined. Finally, again by conjecture \ref{KMconj}
 they depend only on the numbers
$l,k,r,b,w$.
\end{defn}

\begin{lem}
For all $(l,k,r,b,w)$ with $b\ge 0$, if the left hand side of the
equations below is well-defined, then  the right hand side is also, and
\begin{align*}
\tag*{$(0)_s$} P(l,k,r,b,w)&=\sum_{n\in \Z}P(l,k,r,b+1,w-n^2),\\
\tag*{$(1)_s$} P(l,k,r,b,w)&
=\sum_{n\in \Z}(-1)^n(n+1/2)P(l+1,k,r,b+1,w-(n+1/2)^2),\\
\tag*{$(2)_s$} \sum_{n\in \Z} n^2 P(l,k,r,&b+1,w-n^2)=2\sum_{n\in \Z}
P(l-2,k+1,r,b+1,w-n^2),\\
\tag*{$(3)_s$} P(l,k,r+1,b,w)&
=\sum_{n\in \Z}(-1)^{n+1}\Big((n+1/2)^3P(l+3,k,r,b+1,w-(n+1/2)^2)\\
& \quad -6(n+1/2)P(l+1,k+1,r,b+1,w-(n+1/2)^2)\Big).
\end{align*}
\end{lem}
\begin{pf}
Take
$(l,k,r,b,w)$ such that there exists a class $\xi\in H^2(X(b),\Z)$
with $w:=\xi^2/4<0$.
Let $N:=l+2k+2r$. We can assume  that $\xi$ defines a wall of type $(N)$
(otherwise both sides of $(0)_s$--$(3)_s$
are trivially zero).

Assume first that $b_2(X)>2$ and that in addition the intersection
form on $H_2(X,\Z)$ is odd, or  $b>0$.
Then we can find an $\eta$ which is not divisible in $H^2(X(b),\Z)$ with
$\eta^2=\xi^2$.
(The intersection form is $(1)\oplus(-1)^{\oplus b_2(X(b))-1}$,
therefore we can find orthogonal classes $h$, $e_1$, $e_2$ with
$Q_{X(b)}(h)=1=-Q_{X(b)}(e_1)=-Q_{X(b)}(e_2)$, and we put
$\eta:=nh+(n+1)e_1$ (resp. $\eta:=nh+(n+1)e_1+e_2$)
if $\xi^2=-(2n+1)$ (resp. $\xi^2=-(2n+2)$) for $n\in\Z_{\ge 0}$.)

We can therefore assume that $\xi$ is not divisible in $H^2(X(b),\Z)$.
Let $\cc_-$ and $\cc_+$ be the two chambers separated
by $W^\xi$, with $\xi\cdot \check a_-<0<\xi\cdot \check a_+$
for $a_-\in \cc_-$ and $a_+\in \cc_+$.
Assume that
$N+3+4\xi^2<0$. Then $W^{X(b)}_{w,N}(a_-,a_+)=\{\xi\}$.
Therefore we can replace
$\Phi^{X(b),\cc_+}_{c_1,N}-\Phi^{X(b),\cc_-}_{c_1,N}$ in $(0)_r$--$(3)_r$ by
$(-1)^{\varepsilon(c_1,\xi,N)}
\delta^{X(b)}_{\xi,N}$.
Now we apply lemma \ref{combin} and the definition of the
$P(l,k,r,b,w)$ to obtain the result.

Finally, if $m:=N+3+4\xi^2\ge 0$ we make induction over $m$.
So we assume that the result is true for all $m'< m$.
Then $W^{X(b)}_{w,N}(a_-,a_+)=\{\xi\}\cup W_{m}$, where
the classes $\eta\in W_{m}$ satisfy $N+3+\eta^2< m$.
So by induction the result holds for all $\eta\in W_{m}$ and thus by
lemma \ref{wallchange1} also for $\xi$.
\end{pf}

We want to use the $P(l,k,r,b,w)$ as the coefficients of a
power series, which should solve a system of  differential equations.
This  does not work directly,
 because in the moment we have only coefficients
with $w<0$.
So we have to "complete" the coefficients, i.e.
to  define the $P(l,k,r,b,w)$
for all $l,k,r,b,w$ by making use of relation $(1)_s$.

\begin{defn}\label{series}
For all $l,k,r,b\in \Z$ and all  $w\in {1\over 4}\Z$
define inductively $P(l,k,r,b,w)$ by
\begin{enumerate}
\item
If $w=\xi^2/4<0$ for $\xi\in H^2(X(b),\Z)$, then apply definition
 \ref{pcoeff}.
\item We put
$$P(l,k,r,b,w):=\sum_{n\in\Z}(-1)^n(n+1/2)P(l+1,k,r,b+1,w-(n+1/2)^2),$$
whenever the right hand side is already defined inductivly by (1) and (2).
Note that the sum is again finite.
\end{enumerate}
We check that the $P(l,k,r,b,w)$ are well-defined.
For this we have to see (a), that (1) and (2) give the same $P(l,k,r,b,w)$
whenever both apply, but this is the contents of relation
$(1)_s$; and (b),  that the above definition determines
$P(l,k,r,b,w)$ for each $5$-tuple $(l,k,r,b,w)\in \Z^4\times {1\over 4}\Z$.
If $w\le 0$, then there exist on $X(1)$ for all $n\in \Z$ classes $\eta_n$
with $\eta_n^2=4w-(2n+1)^2<0$ (as the intersection form on $X(1)$
is odd and of rank $\ge 3$), and thus $P(l,k,r,b,w)$ is defined by  (2).
Now assume that $P(l,k,r,b,w')$ is defined for all $l,k,r,b$ and all
$w'<w$. Then we use again  (2) to define $P(l,k,r,b,w)$.
We put
$$\Lambda_X(L,Q,x,t,\tau):=
\sum_{(l,k,r,b)\in \Z^4}\sum_{w\in {1\over 4}\Z} P(l,k,r,b,w){
L^lQ^kx^rt^bq^w\over l!k!r!b!},$$
where again $\tau\in \Bbb H$ and $q=e^{2\pi i\tau}$.
\end{defn}

$\Lambda_X$ now encodes all the wall-crossing formulas
for all blowups of $X$.

\begin{rem}
Let $\xi\in H^2(X(b),\Z)$ be a class with $\xi^2<0$.
Then for all $\alpha\in  H_2(X(b),\Q)$
$$\delta_{\xi}^{X(b)}(\exp(\alpha z+p x))=\hbox{res}_{q=0}
\left({\partial^k \over \partial t^k}(q^{-\xi^2/4}\Lambda_X((\xi/2\cdot\alpha)
z,Q_X(\alpha)z^2,x,t,\tau)dq/q\right)\Big|_{t=0}.$$
\end{rem}
\begin{pf}
This follows directly from the definition.
\end{pf}

\begin{lem}
 $\Lambda_X$ satisfies the differential equations
\begin{align}
 \theta(\tau){\partial\over \partial t}\Lambda_X&=\Lambda_X,\tag*{$(0)_d$}\\
\eta(2\tau)^3 {\partial\over \partial L}{\partial\over \partial t}
\Lambda_X&=\Lambda_X,\tag*{$(1)_d$}\\
 2\theta(\tau){\partial\over \partial Q}\Lambda_X&=(q{d\over
dq}\theta(\tau)){\partial^2\over \partial L^2}\Lambda_X,\tag*{$(2)_d$}\\
{\partial\over \partial x}\Lambda_X&=(q{d\over dq}\eta(2\tau)^3)
{\partial^3\over \partial L^3}{\partial\over \partial t}\Lambda_X-
6\eta(2\tau)^3
{\partial\over \partial L}
{\partial\over \partial Q}{\partial\over \partial t}\Lambda_X.\tag*{$(3)_d$}
\end{align}
\end{lem}

\begin{pf}
We first want to see that the relations
$(0)_s$--$(3)_s$ hold for all $(l,k,r,b,w)\in \Z^4\times {1\over 4}\Z$,
i.e. that the recursive definition
 is compatible with $(0)_s$--$(3)_s$. The proof is similar
in all cases, so we just do $(0)_s$.
We assume that $(0)_s$ holds for all $(l,k,r,b,w')$ with $w'<w$.
Then we get
\begin{eqnarray*}
P(l,k,r,b,w)&=&\sum_{n\in \Z}(-1)^n(n+1/2)P(l+1,k,r,b+1,w-(n+1/2)^2)\\
&=&\sum_{n,m\in \Z}(-1)^n(n+1/2)P(l+1,k,r,b+2,w-(n+1/2)^2-m^2)\\
&=&\sum_{m\in \Z}P(l+1,k,r,b+2,w-m^2).
\end{eqnarray*}
We now  translate
$(0)_s$--$(3)_s$ into differential equations
$(0)_d$--$(3)_d$:
\begin{eqnarray*}
\Lambda_X&=&\sum_{(l,k,r,b,w)} P(l,k,r,b,w){
L^lQ^kx^rt^bq^w\over l!k!r!b!}\\
&=&\sum_{(l,k,r,b,w)} \sum_{n\in\Z}
P(l,k,r,b,w){L^lQ^kx^rt^{b-1}q^{w+n^2}\over l!k!r!(b-1)!}\\
&=&\theta(\tau){\partial\over \partial t}\Lambda_X.
\end{eqnarray*}
Similarly we get
\begin{eqnarray*}
\Lambda_X&=&
\sum_{(l,k,r,b,w)} \sum_{n\in\Z}
(-1)^n(n+1/2)P(l,k,r,b,w){L^{l-1}Q^kx^rt^{b-1}q^{w+(n+1/2)^2}\over
(l-1)!k!r!(b-1)!}\\
&=&\sum_{n\in\Z}(-1)^n(n+1/2)q^{(n+1/2)^2}
{\partial\over \partial L}{\partial\over \partial t}\Lambda_X,
\end{eqnarray*}
and $(1)_d$ follows from remark \ref{modular}. Furthermore
\begin{eqnarray*}
0&=&\sum_{(l,k,r,b,w)} P(l,k,r,b,w)\sum_{n\in\Z}
\left({n^2L^{l-2}Q^k\over (l-2)!k!}-{2L^{l}Q^{k-1}\over
l!(k-1)!}\right){x^rt^{b-1}q^{w+n^2}\over r!(b-1)!}\\
&=&(q {d\over dq}\theta(\tau))
{\partial^2\over \partial L^2}{\partial\over \partial t}\Lambda_X
-2{\partial\over \partial Q}{\partial\over \partial t}\Lambda_X.
\end{eqnarray*}
Finally we get
\begin{eqnarray*}
{\partial\over \partial x}\Lambda_X
&=&\sum_{(l,k,r,b,w)}P(l,k,r+1,b,w){L^{l}Q^kx^rt^{b}q^{w}\over l!k!r!b!}\\
&=&\sum_{(l,k,r,b,w)}P(l,k,r,b,w) \sum_{n\in\Z}
(-1)^{n+1}\left({(n+1/2)^3L^{l-3}Q^k\over (l-3)!k!}\right.\\
&&\qquad
-\left.{6(n+1/2)L^{l-1}Q^{k-1}\over
(l-1)!(k-1)!}\right){x^rt^{b-1}q^{w+(n+1/2)^2}\over r!(b-1)!}\\
&=&(q{d\over dq}\eta(2\tau)^3)
{\partial^3\over \partial L^3}{\partial\over \partial t}\Lambda_X-
6\eta(2\tau)^3{\partial\over \partial L}{\partial\over \partial
Q}{\partial\over \partial t}\Lambda_X.
\end{eqnarray*}
\end{pf}

\begin{lem}\label{expex}
Putting $\lambda_X(\tau):=\Lambda_X(0,0,0,0,\tau)$  we obtain
$$\Lambda_X=\exp\Big(L/f(\tau)-(Q/2)(2G_2(2\tau)
+e_3(2\tau))/f(\tau)^2-3xe_3(2\tau)/f(\tau)^2-
\theta(\tau)t\big)\lambda_X(\tau).$$
\end{lem}
\begin{pf}
Using remark \ref{modular} we can reformulate
$(0)_d$--$(3)_d$ as
\begin{eqnarray*}
{\partial\over \partial t}\Lambda_X&=&\Lambda_X/\theta(\tau),\qquad
{\partial\over \partial L}\Lambda_X=\Lambda_X/f(\tau),\\
{\partial\over \partial Q}\Lambda_X&=&{1\over
2}q\,(d\log_q(\theta(\tau)))\Lambda_X/f(\tau)^2,\\
{\partial\over \partial x}\Lambda_X
&=&q\,(d\log_q(\eta(2\tau)^3))\Lambda_X/f(\tau)^2
-3q(d\log_q(\theta(\tau)))\Lambda_X/f(\tau)^2.
\end{eqnarray*}
So the result follows by remark \ref{modular}.
\end{pf}

To finish the proof of  theorem \ref{mainthm} we now only have to identify
$\lambda_X$.

\begin{lem}\label{norm}
$\lambda_X=\Delta(2\tau)\exp(-\theta(\tau)\sigma(X))/f(\tau)^{11}=
f(\tau)\Delta(2\tau)^2/(\Delta(\tau)\Delta(4\tau))
\,\exp(-\theta(\tau)\sigma(X)).$
\end{lem}
\begin{pf}
We first show that is is enough to prove this result  in case $X=\P_1\times
\P_1$.  We note that by lemma \ref{expex} the statement for a variety $Y$ and
$\widehat Y=Y\#\overline \P_2$ are equivalent. It is therefore enough to show
it for $\widehat X$.
$\widehat X$ has odd intersection form and $a:=b_2(\widehat X)-1\ge 2$.
So it is homotopy-equivalent to
$\P_2\# a\overline \P_2=(\P_1\times\P_1)\#(a-1)\overline\P_2$.
As $\delta_{\xi}^X$ only depends on the homotopy type of $X$,
it is enough to show
the result for $\P_1\times\P_1$. This is in fact the only time in our
argument where we use that $\delta_{\xi}^X$ depends on the homotopy type
of $X$, rather then on $X$ itself.

Let $F,G\in H^2(\P_1\times\P_1,\Z)$ be the classes of the fibres
of the two projections to $\P_1$.
Let $k\in \Z_{>0}$ and $N:=4k-1$.
Then lemma \ref{qin} gives that  $\Phi^{\P_1\times\P_1,F+\epsilon
G}_{F+G,N}=\Phi^{\P_1\times\P_1,G+\epsilon F}_{F+G,N}=0$
for all sufficiently small $\epsilon>0$.
In particular we have for all $k>0$
$$(-1)^{k+1}\sum_{\xi\in
W^{\P_1\times\P_1}_{f+g}(F,G)}(-1)^{\varepsilon(F+G,\xi,4k-1)}
\delta_\xi^{\P_1\times\P_1}(2\check G^{4k-1})=0.$$
Here $W^{\P_1\times\P_1}_{f+g}(F,G)=\big\{(2n-1)F-(2m-1)G\bigm| n,m\in
\Z_{>0}\big\}$,
and $(-1)^{k+1+\varepsilon(F+G,(2n-1)F-(2m-1)G,4k-1)}=(-1)^{n+m}$.
Applying lemma \ref{expex} we get
\begin{equation}\hbox{res}_{q=0}\Big(\sum_{n,m\in\Z_{>0}}(-1)^{n+m}q^{{1\over
2}(2n-1)(2m-1)}
(2n-1)^{4k-1}f(\tau)^{-4k+1}\lambda_{\P_1\times\P_1}(\tau){dq\over q}\Big)
=0.\tag{*}
\end{equation}
Note that, by remark \ref{power}, $\lambda_{\P_1\times\P_1}=
q^{-3/4}\bar\lambda$ where $\bar\lambda=\sum l_i q^i$ is a power series in $q$.
Also $f(\tau)=q^{1/4}\bar f$ with $\bar f$ a  power series
in $q$ with constant term $1$.
It is well-known (\cite{Ko},\cite{K-M}) that
$\delta^{\P_1\times\P_1}_{F-3G}((2\check G)^3)=1$. Thus we get $l_0=1$.
$(*)$ gives for each $k\ge 1$ the recursive relation
$$\sum_{n,m>0} (-1)^{n+m}(2n-1)^{4k-1}\hbox{Coeff}_{q^{k-2nm+n+m}}(
\bar \lambda/\bar f^{4k-1})=0,$$
i.e., putting $\lambda_k:=\sum_{j<k} l_jq^j$ we obtain
$$l_k=-\sum_{n,m>0} (-1)^{n+m}(2n-1)^{4k-1}\hbox{Coeff}_{q^{k-2nm+n+m}}(\bar
\lambda_k/\bar f^{4k-1}).$$
So we see that $\lambda_{\P_1\times\P_1}$ is uniquely determined by $(*)$.
We put
$$H_k(\tau):=\sum_{n,m\in\Z_{>0}}(-1)^{n+m}q^{{1\over 4}(2n-1)(2m-1)}
(2n-1)^{4k-1}\Delta(\tau)/f(\tau/2)^{4k+10}.$$
Then the lemma follows from the following
lemma (the  proof of which
is due to Don Zagier).
\end{pf}

\begin{lem}\label{residue} $\hbox{res}_{q=0}H_k(\tau){dq\over q}=0$.
\end{lem}
\begin{pf}
We start by rewriting $H_k(\tau)$.
\begin{eqnarray*}
\sum_{n,m>0}(-1)^{n+m}q^{(n-1/2)(m-1/2)}(2n-1)^{4k-1}
=\sum_{d\ \text{\rm odd}}^\infty(-1)^{(d-1)/2}\sigma_{4k-1}(d)q^{d/4}\\
={1\over 2i}(G_{4k}((\tau+1)/ 4)-G_{4k}((\tau-1)/ 4))=:\widetilde G_{4k}(\tau),
\end{eqnarray*}
where $G_{4k}(\tau)$ is the Eisenstein series.
We write
$\phi:=f(\tau/2)^{2}=\left(\eta(\tau/2)\eta(2\tau)/\eta(\tau)\right)^4.$
Then we have $H_k(\tau)={\widetilde G_{4k}(\tau)\Delta(\tau)/\phi^{2k+5}}$.
We want to show that $H_k(\tau)$ is a modular form of weight $2$
for the $\theta$-group
$$\Gamma_{\theta}:=\big\{ A\in SL(2,\Z)\bigm| A\equiv \left(\begin{matrix}
1& 0\\0 & 1\end{matrix}\right) or  A\equiv \left(\hbox{$\begin{matrix}
0 & 1\\1& 0\end{matrix}$}\right)\hbox{ modulo } 2 \big\}.$$
The operation of $\Gamma_{\theta}$ is generated by
$\tau\mapsto \tau+2$ and $\tau\mapsto {-1/\tau}$.
We see that $\widetilde G_{4k}(\tau+2)=
-\widetilde G_{4k}(\tau)$. Now we write
$${(-1/\tau+1)/4}={(\tau-1)/4 \over 4(\tau-1)/4+1}, \ \
{(-1/\tau-1)/4}={(\tau+1)/4 \over -4(\tau+1)/4+1},$$
and use  that
$G_{4k}(\tau)$ is a modular form of weight $4k$ for $SL(2,\Z)$,
to obtain that
$$\widetilde G_{4k}(-1/\tau)={1\over 2i}\big(\tau^4G_{4k}((\tau-1)/ 4)
-\tau^4 G_{4k}((\tau+1)/ 4)\big)=-\tau^4\widetilde G_{4k}(\tau).
$$
Furthermore we see  by  $\phi(\tau)^6=\Delta(\tau/2)\Delta(2\tau)/
\Delta(\tau)$, that $\phi(-{1/\tau})^6=
\tau^{12 }\phi(\tau)^6$, i.e. $\phi(-{1/\tau})=\omega
\tau^{2 }\phi(\tau)$  for a $6$-th root of unity $\omega$.
Putting $\tau:=i$ we get $\phi(i)=-\omega\phi(-{1/ i})$, i.e. $\omega=-1$.
We also obviously have $\phi(\tau+2)=-\phi(\tau)$.
Putting this together and using the fact that $\Delta(\tau)$
is a modular form  of weight $12$ for $SL(2,\Z)$ we finally
see that
$H_k(\tau)$ is a modular form of weight $2$ for $\Gamma_\theta$.
In other words
$H_k(\tau){dq/ q}=2\pi i H_k(\tau){d\tau }$ is a differential form on
the rational curve $\overline {\Bbb H/\Gamma_\theta}$,  holomorphic out of the
cusps
$\tau=1$ and $\tau=\infty$ (i.e. $q=0$).
We show that $H_k(\tau)$ is holomorphic at $\tau=1$.
$\Delta(\tau)$ and $\widetilde G_{4k}(\tau)$ are obviously holomorphic
at $\tau=1$. We now put
 $\tau:=1-1/z$ and use again  that $\Delta(\tau)$ is a modular form of weight
$12$ for $SL(2,\Z)$ and write
$1/ 2-1/ (2 z)={(z-1)/2\over 2(z-1)/2+1}$ to obtain
$$\phi(\tau)^6={\Delta(1/ 2-
1/ (2 z))\Delta(-2/z)\over
\Delta(-1/z)}=
{z^{12}\Delta((z-1)/ 2)(z/ 2)^{12}\Delta(z/ 2)\over
z^{12}\Delta(z)}=-(z/ 2)^{12}{\Delta(z)^2\over \Delta(2z)}.
$$
So for $z=\infty$, (i.e. $\tau=1$) the modular form
$\phi(\tau)$ is holomorphic and does not vanish. Thus also
$H_k(\tau)$ is holomorphic at $\tau=1$.
Thus  the residue theorem implies that
$\hbox{res}_{q=0}( 2\pi iH_k(\tau){d\tau  })=0.$
\end{pf}

\begin{rem} \label{homotop}
As noted above, we have used that by  conjecture \ref{KMconj}
$\delta_{\xi,N}^X$ depends only on the homotopy type $X$ rather
then just on $X$ only in the reduction above to $\P_1\times\P_1$.
In particular, without assuming this, our proof still shows
theorem \ref{mainthm} for $X$ a rational surface, and therefore also
theorem \ref{donp2}.
\end{rem}

\section{Possible generalizations}
It should be possible to prove the blowup formulas and  also
conjecture \ref{KMconj}  for $4$-manifolds $X$ with
$b_+(X)=1$ and $b_1(X)=0$ (i.e dropping the assumption that $X$ is
simply-connected).
If we assume these generalizations, then all our arguments
in the proof of theorem \ref{mainthm} work in this more general case
except for the reduction to $\P_1\times \P_1$
at the beginning of the proof of lemma \ref{norm}.
So we get

\begin{cor}
Assume that the blowup formulas \ref{blowup}
and  conjecture \ref{KMconj} hold for all for $4$-manifolds $Y$
with $b_+(Y)=1$ and $b_1(Y)=0$.
Then for all  $X$ with $b_+(X)=1$ and $b_1(X)=0$,
all $\xi$ in $H^2(X,\Z)$ with $\xi^2<0$ and all
$\alpha\in H_2(X,\Q)$ we have
$$\delta_{\xi}^{X}(\exp(\alpha z+p x))=\hbox{res}_{q=0}(g^X_{\xi}(\alpha z,
x,\tau)\lambda_{[X]}(\tau)\Delta(\tau)\Delta(4\tau)/(f(\tau)\Delta(2\tau)^2))
dq/q,$$
where $g^X_{\xi}$ is the generating function from theorem \ref{mainthm}
and  $\lambda_{[X]}(\tau)$ is
$q^{-3/4}$ multiplied with an unknown power series $\bar\lambda_{[X]}(q)$
in $q$, which depends only on the equivalence class $[X]$,
where $X$ and $Y$ are equivalent if
$X\#k\overline\P_2$ and $Y\#k\overline\P_2$ are homotopy equivalent
for some $k$.
\end{cor}

The results of \cite{E-G1} suggest that the dependence of
$\lambda_{[X]}(\tau)$ on $X$ should be very simple.

\begin{conj}
$\lambda_{[X]}(q)=n_2f(\tau)\Delta(2\tau)^2/(\Delta(\tau)\Delta(4\tau)),$
for $n_2$  the number of $2$-torsion points in $H^2(X,\Z)$.
\end{conj}

\end{document}